\documentclass{amsart}

\DeclareMathSymbol{\twoheadrightarrow}  {\mathrel}{AMSa}{"10}

\def\GG{{\mathcal G}}

\def\Q{{\mathbb Q}}
\def\Z{{\mathbb Z}}
\def\C{{\mathbb C}}

\def\F{{\mathbb F}}

\def\Sn{{\mathbf S}_n}
\def\An{{\mathbf A}_n}
\def\RR{{\mathfrak R}}

\def\Perm{\mathrm{Perm}}
\def\Gal{\mathrm{Gal}}

\def\End{\mathrm{End}}
\def\Aut{\mathrm{Aut}}
\def\Hom{\mathrm{Hom}}

       \def\Lie{\mathrm{Lie}}

\def\I{{\mathcal I}}

\def\J{{\mathcal J}}

\def\ST{{\mathbf S}}

\def\fchar{\mathrm{char}}

\def\SL{\mathrm{SL}}
\def\Sp{\mathrm{Sp}}
    \def\SA{\mathrm{SA}}
     \def\SPA{\mathrm{ASP}}
\def\pr{\mathrm{pr}}
\def\M{\mathrm{M}}

                      \def\G{{\mathcal G}}

\def\dim{\mathrm{dim}}

\def\OC{{\mathcal O}}

\def\a{{\mathfrak a}}
\def\b{{\mathfrak b}}

\newtheorem{thm}{Theorem}[section]
\newtheorem{lem}[thm]{Lemma}
\newtheorem{cor}[thm]{Corollary}

\theoremstyle{definition}
\newtheorem{defn}[thm]{Definition}
\newtheorem{ex}[thm]{Example}
\newtheorem{rem}[thm]{Remark}

\title[Endomorphisms of superelliptic jacobians]
{Endomorphisms of superelliptic jacobians}
\author[Yuri\ G.\ Zarhin]{Yuri\ G.\ Zarhin}
\address{Department of Mathematics, Pennsylvania State University,
University Park, PA 16802, USA}
\address{Institute for Mathematical Problems in Biology, Russian Academy of
Sciences, Pushchino, Moscow Region, Russia}
 \email{zarhin\char`\@math.psu.edu}

\begin{document}

\begin{abstract}
  Let $K$ be a field of characteristic zero, $n\ge 5$ an integer, $f(x)$ an irreducible
polynomial over $K$ of degree $n$, whose Galois group contains a doubly
transitive simple non-abelian group. Let $p$ be an odd prime,
$\Z[\zeta_p]$ the ring of integers in the $p$th cyclotomic field,
 $C_{f,p}:y^p=f(x)$  the  corresponding superelliptic curve and
$J(C_{f,p})$ its jacobian. Assuming that either $n=p+1$ or $p$
does not divide $n(n-1)$, we prove that the ring of all
endomorphisms of $J(C_{f,p})$ coincides with $\Z[\zeta_p]$.
The same is true if $n=4$, the Galois group of $f(x)$ is the full symmetric group $\ST_4$
 and $K$ contains a primitive $p$th root of unity.

2000 Math. Subj. Class: Primary 14H40; Secondary 14K05, 11G30,
11G10

\end{abstract}

\maketitle

This is a corrected version of \cite{ZarhinMZ2}.  Inaccuracies in the statements of Theorems 1.1(ii), 3.12(ii), 5.2(ii) and Remarks 1.3,  3.2 of  \cite{ZarhinMZ2} and a gap in the proof of Theorem 3.12(ii) \cite[p. 702]{ZarhinMZ2} (caused by improper use of Theorem 4.3.2 of \cite{Herstein} in \cite[p. 697]{ZarhinMZ2}) are fixed. Several typos have been corrected. The reference to  \cite[Th. 4.3.2]{Herstein} is moved to the end of Section \ref{endom} (Lemma \ref{herst}).

 I am deeply grateful to Jiangwei Xue, who pointed out the gap.

\section{Definitions and statements}
\label{one}
\label{endo} Throughout this paper $K$ is a field and $K_a$  its
algebraic closure. We write $\Gal(K)$ for the absolute Galois
group $\Gal(K_a/K)$ and $\ell$ for a prime different from
$\fchar(K)$. If $X$ is an abelian variety of positive dimension
 over $K_a$ then $\End(X)$ stands for the ring of all its
$K_a$-endomorphisms and $\End^0(X)$ for the corresponding
$\Q$-algebra $\End(X)\otimes\Q$. For every  abelian variety over
$Y$ over $K_a$  we write $\Hom(X,Y)$ for the group of all
$K_a$-homomorphisms from $X$ to $Y$.

Let $p$ be a prime and $\zeta_p \in \C$ be a primitive $p$th root of unity,
$\Q(\zeta_p)\subset \C$ the
$p$th cyclotomic field and $\Z[\zeta_p]$ the ring of integers in $\Q(\zeta_p)$.

Let us assume that $\fchar(K)\ne p$ and $K$ contains a primitive
$p$th root of unity $\zeta$. Let $f(x) \in K[x]$ be a polynomial
of degree $n\ge 3$ without multiple roots,
 $\RR_f\subset K_a$ the ($n$-element) set of roots of $f$ and $K(\RR_f)\subset K_a$
the splitting field of $f$. We write $\Gal(f)=\Gal(f/K)$ for the
Galois group $\Gal(K(\RR_f)/K)$ of $f$; it permutes the roots of
$f$ and may be viewed as a certain permutation group of $\RR_f$,
i.e., as  a subgroup of the group $\Perm(\RR_f)\cong\Sn$ of
permutations of $\RR_f$.

 Let $C_{f,p}$ be a smooth projective model of the smooth affine $K$-curve
$y^p=f(x)$; its genus is equal to $(n-1)(p-1)/2$ if $p$ does {\sl not} divide $n$. The map
 $(x,y) \mapsto (x, \zeta y)$
gives rise to a non-trivial birational $K$-automorphism $\delta_p:
C_{f,p} \to C_{f,p}$ of period $p$. The jacobian $J(C_{f,p})$ of
$C_{f,p}$ is an abelian variety that is defined over $K$.  By Albanese
functoriality, $\delta_p$ induces an automorphism of $J(C_{f,p})$
which we still denote by $\delta_p$. It is known \cite[p.~149]{Poonen}, \cite[p.~458]{SPoonen})
(see also \cite{ZarhinM,ZarhinPisa})
that $\delta_p$ satisfies the $p$th cyclotomic equation, i.e.,
${\delta_p}^{p-1}+\cdots +\delta_p+1=0$
 in $\End(J(C_{f,p}))$. This gives rise to the embeddings
 $$\Z[\zeta_p] \hookrightarrow \End_K(J(C_{f,p}))\subset \End(J(C_{f,p})),\
 \zeta_p\mapsto \delta_p,$$
$$\Q(\zeta_p) \hookrightarrow \End^0_K(J(C_{f,p}))\subset \End^0(J(C_{f,p})),\
 \zeta_p\mapsto \delta_p,$$
that send $1$ to the identity automorphism of $J(C_{f,p})$.
In particular,  if $\Q[\delta_p]$ is the
$\Q$-subalgebra of $\End^0(J(C_{f,p}))$ generated by $\delta_p$
then the latter embedding establishes a canonical isomorphism
$$\Q(\zeta_p)\cong \Q[\delta_p]$$
that sends $\zeta_p$ to $\delta_p$. In a series of papers
\cite{ZarhinMRL,ZarhinMMJ,ZarhinBSMF,ZarhinCrelle,ZarhinCamb,ZarhinSb,ZarhinM,ZarhinPisa},
the author discussed the structure of $\End^0(J(C_{f,p}))$.
 In particular, he proved that if $\fchar(K)=0, \ n \ge
5$ and $\Gal(f)$ coincides either with  full symmetric group $\Sn$
or with  alternating group $\An$
then  $\End(J(C_{f,p}))=\Z[\zeta_p]$. (See also \cite{ZarhinTexel,ZarhinMRL2,ZarhinPAMS,ZarhinSh,ZarhinMZ}.)

The aim of this  paper is to deal with arbitrary doubly transitive Galois groups for almost all $p$ (if  $n$ is given).
 For example, we extend the result about $\End(J(C_{f,p}))$ to the case when $n=4, \Gal(f)=\ST_4$,
  a prime $p$ is odd and $K$ contains a primitive $p$th root of unity (see below).

The following statement is the main result of this paper.

 \begin{thm}
 \label{main}
Suppose that $K$ has characteristic zero, $n\ge 4$ and $p$ is an
odd prime that does not divide $n$. Suppose that $K$ contains a
primitive $p$th root of unity and
  $\Gal(f)$ contains a doubly transitive
subgroup $\G$ that enjoys the following property: if $p$ divides $n-1$ and $n\ne p+1$ then $\G$
 does not contain a proper subgroup, whose index divides $(n-1)/p$.

 Then:

 \begin{itemize}
\item[(i)]
$\End^0(J(C_{f,p}))$ is a simple $\Q$-algebra, i.e., $J(C_{f,p})$ is isogenous over $K_a$
to a self-product of an absolutely simple abelian variety.
\item[(ii)]

Let us assume  that  $\G$ does not contain a proper normal subgroup, whose  index divides $n-1$. Assume also that  either $n=p+1$ or  $p$ does not divide $n-1$.
Then
 $$\End^0(J(C_{f,p}))=\Q[\delta_p]\cong\Q(\zeta_p), \ \End(J(C_{f,p}))=\Z[\delta_p]\cong\Z[\zeta_p].$$
 \end{itemize}
\end{thm}

\begin{rem}
See \cite{DM,Mortimer} for the list of known doubly transitive permutation groups.
\end{rem}

\begin{rem}
\label{proper}
Since $\G$ is doubly transitive, its order is divisible by $n(n-1)$. In particular,
every subgroup in $\G$, whose index divides $n-1$ is either nontrivial proper or $\G$ itself.
\end{rem}

\begin{rem}
The case of $n=3$ is discussed in \cite{ZarhinPisa}; see also
\cite{Ribet3}. (Notice that the only doubly transitive subgroup of
$\mathbf{S}_3$ is $\mathbf{S}_3$ itself, which contains a normal
subgroup of index $2=3-1$.)
\end{rem}

Let us consider the case of $n=4$.
 The only doubly transitive subgroups of $\mathbf{S}_4$ are $\mathbf{S}_4$ and  $\mathbf{A}_4$:
 they both contain a  subgroup of index $3=4-1$. The only  subgroup in $\mathbf{A}_4$ of index
 $3$ is normal. However, none of the index $3$ subgroups in $\mathbf{S}_4$ is normal.
(Indeed, an index $3$ normal subgroup of $\mathbf{S}_4$ is,
obviously, a normal Sylow $2$-subgroup and therefore must contain
all the transpositions, which could not be the case.)
  Now Theorem \ref{main}
  implies the following statement.

  \begin{cor}
 \label{main4}
Suppose that $K$ has characteristic zero, $f(x)\in K[x]$ is an irreducible
 quartic polynomial, whose Galois group  is the full symmetric group $\mathbf{S}_4$.
  If $p$ is an odd prime  and $K$ contains a primitive $p$th root of unity then
$$\End^0(J(C_{f,p}))=\Q[\delta_p]\cong\Q(\zeta_p), \ \End(J(C_{f,p}))=\Z[\delta_p]\cong\Z[\zeta_p].$$
\end{cor}

\begin{ex}
\label{quartic}
 Let $p=3$ and $f(x)=\sum_{i=0}^4 f_i x^i \in K[x]$
be an irreducible degree $4$ polynomial over $K$ with Galois group
$\mathbf{S}_4$. Then $f_0\ne 0, f_4\ne 0$ and the plane projective
quartic
$$y^3z-\sum_{i=0}^4 f_i x^i z^{4-i}=0$$
is a smooth projective model of $C_{f,3}$, whose only {\sl
infinite} point is $(0:1:0)$. In particular, $C_{f,3}$ is {\sl
not} hyperelliptic. It follows from Corollary \ref{main4} that
$$\End(J(C_{f,3}))=\Z\left[\frac{-1+\sqrt{-3}}{2}\right]$$
if $\fchar(K)=0$ and $K$ contains a primitive cubic root of unity.

For example, the polynomial $f(x)=x^4-x-1$ has Galois group
$\mathbf{S}_4$ over $\Q$ \cite[p. 42]{Serre}. Let $L$ be the
splitting field of $f(x)$ over $\Q$, so that $L/\Q$ is a Galois
extension with $\Gal(L/\Q)=\mathbf{S}_4$. Since $\mathbf{S}_4$ has
exactly one subgroup of index $2$ (namely, $\mathbf{A}_4$), the
field $L$ contains exactly one quadratic (sub)field, namely,
$\Q(\sqrt{\Delta})$ where $\Delta$ is the discriminant of $f(x)$.
Recall that the discriminant of a quartic polynomial $x^4+ax+b$ is
$-27 a^4+256 b^3$ \cite[Sect. 64]{V}. It follows that
$\Delta=-283$ and $L$ contains $\Q(\sqrt{-283})$. This implies
that $L$ does {\sl not} contain $K:=\Q(\sqrt{-3})$ and therefore
$L$ and $K$ are linearly disjoint over $\Q$. Hence the Galois
group of $f(x)$ over $K$ is still $\mathbf{S}_4$. Since $K$
contains a primitive cubic root of unity, we conclude that if $J$
is the jacobian of the smooth plane projective quartic
$$y^3z-x^4+xz^3+z^4=0$$ then
$\End(J)=\Z\left[\frac{-1+\sqrt{-3}}{2}\right]$.

 Notice that one
may use  degree seven polynomials with Galois group $\mathbf{S}_7$
or $\mathbf{A}_7$ in order to construct smooth plane quartics,
whose jacobians have no nontrivial endomorphisms \cite{ZarhinVE}.
\end{ex}

Recall that one may realize $\mathbf{S}_4$ as the group of affine
transformations of the affine plane over $\F_2$. This suggests the
following generalization of Corollary \ref{main4}.

\begin{thm}
\label{mainA} Suppose that there exists a prime $\ell\ne p$ and
positive integers $r$ and $d$ such that $d>1$ and $n=\ell^{rd}$.
Let $\F$ be the finite field $\F_{\ell^r}$. Let $\F^d$ be the
$d$-dimensional coordinate $\F$-vector space and $\SA(d,\F)$ be
the group of all special affine transformation
$$v \mapsto A(v)+t, \ t\in \F^d, A\in SL(d,\F)$$
of $\F^d$. Suppose that one may identify $\RR_f$ and $\F^d$ in
such a way that $\Gal(f)$  becomes a permutation group of $\F^d$
that contains $\SA(d,\F)$.

Suppose that $\fchar(K)=0$  and $K$ contains a primitive $p$th
root of unity. If either $n=p+1$ or $p$ does not divide $n-1$ then
$$\End^0(J(C_{f,p}))=\Q[\delta_p]\cong\Q(\zeta_p), \ \End(J(C_{f,p}))=\Z[\delta_p]\cong\Z[\zeta_p].$$
\end{thm}

\begin{proof}
The permutation group $\G=\SA(d,\F)$ is doubly transitive, because
$\SL(d,\F)$ is transitive on the set of all non-zero vectors in
$\F^d$ (recall that $d>1$). The group $T\cong \F^d$ of all
translations of $\F^d$ is normal in $\SA(d,\F)$ and
$\G/T=\SL(d,\F)$.  Suppose that $\G$ contains a proper normal
subgroup say, $H$, whose index divides $n-1=\ell^{rd}-1$. It
follows that $H$ contains all elements, whose order is a power of
$\ell$. This implies that $H$ contains $T$ and therefore $H/T$ is
a normal subgroup in $\SL(d,\F)$ that contains all elements of
$\SL(d,\F)$, whose order is a power of $\ell$. Since $\SL(d,\F)$
is generated by {\sl elementary matrices} \cite[Ch. XIII, Sect. 8.
Lemma 8.1 and  Sect. 9, Proposition 9.1]{Lang} and every
(non-identity) elementary matrix over $\F$ has multiplicative
order $\ell$, we conclude that $H/T=\SL(d,\F)$ and therefore
$H=\SA(d,\F)$. Since $H$ is proper, we get a contradiction. So
such $H$ does not exist. Now Theorem \ref{mainA} follows from
Theorem \ref{main}.
\end{proof}

Now let us discuss the case of $n\ge 5$ and {\sl simple} $\G$.

\begin{cor}
\label{simple}
Suppose that $\fchar(K)=0$ and $n\ge 5$.  Suppose that $p$ is an odd prime that does not
divide $n$ and such that either $n=p+1$ or $p$ does not divide $n-1$. Suppose that
  $\Gal(f)$ contains a doubly transitive simple non-abelian subgroup  $\G$.
Then
$$\End^0(J(C_{f,p}))=\Q[\delta_p]\cong\Q(\zeta_p), \  \End(J(C_{f,p}))=\Z[\delta_p]\cong\Z[\zeta_p].$$
\end{cor}

\begin{proof}[Proof of Corollary \ref{simple}]
Replacing $K$ by a suitable finite algebraic extension, we may
assume that $\Gal(f)=\G$. Since $\G$ is simple non-abelian,
$\Gal(f)$ would not change if we replace $K$ by an abelian
extension. So, we may assume that $K$ contains a primitive $p$th
root of unity. The simplicity of $\G$ means that it does {\sl not}
have a proper normal subgroup. Now the assertion follows from
Theorem \ref{main}.
\end{proof}

The paper is organized as follows. Section \ref{endom} recalls
elementary but useful results from \cite{ZarhinM} concerning the
endomorphism algebras of abelian varieties. These results are
generalized in Section \ref{abmult} to the case of the centralizer
of a given field that acts on an abelian variety. We pay a special
attention to the action of this field on the differentials of the
first kind (Theorems \ref{multtangent} and \ref{mainAV}).  In
Section \ref{proofm}  we prove Theorem \ref{main}. In Section
\ref{comq} we discuss generalizations of our results to the case
when $p$ is replaced by its power $q$ and the curve involved is
$y^q=f(x)$.

This work was started during the special semester ``Rational and
integral points on higher-dimensional varieties" at the MSRI
(Spring 2006).  I am grateful to the MSRI and the organizers of
this program.

\section{Endomorphism algebras of abelian varieties}
\label{endom}
Let $X$ be an abelian variety of positive dimension over $K$.
 If $n$ is a positive
integer that is not divisible by $\fchar(K)$ then $X_n$ stands for
the kernel of multiplication by $n$ in $X(K_a)$. It is well-known
\cite{MumfordAV} that  $X_n$ is a free $\Z/n\Z$-module of rank
$2\dim(X)$. In particular, if $n=\ell$ is a prime then $X_{\ell}$
is a $2\dim(X)$-dimensional $\F_{\ell}$-vector space.

 If $X$ is defined over $K$ then $X_n$ is a Galois
submodule in $X(K_a)$ and all points of $X_n$ are defined over a
finite separable extension of $K$. We write
$\bar{\rho}_{n,X,K}:\Gal(K)\to \Aut_{\Z/n\Z}(X_n)$ for the
corresponding homomorphism defining the structure of the Galois
module on $X_n$,
$$\tilde{G}_{n,X,K}\subset
\Aut_{\Z/n\Z}(X_{n})$$ for its image $\bar{\rho}_{n,X,K}(\Gal(K))$
and $K(X_n)$ for the field of definition of all points of $X_n$.
Clearly, $K(X_n)$ is a finite Galois extension of $K$ with Galois
group $\Gal(K(X_n)/K)=\tilde{G}_{n,X,K}$. If $n=\ell$  then we get
a natural faithful linear representation
$$\tilde{G}_{\ell,X,K}\subset \Aut_{\F_{\ell}}(X_{\ell})$$
of $\tilde{G}_{\ell,X,K}$ in the $\F_{\ell}$-vector space
$X_{\ell}$.

We write $\End_K(X)$ for the ring of all $K$-endomorphisms of $X$.
We have
$$\Z=\Z\cdot 1_X \subset \End_K(X) \subset \End(X)$$ where $1_X$
is the identity automorphism of $X$. Since $X$ is defined over
$K$, one may associate with every $u \in \End(X)$ and $\sigma \in
\Gal(K)$ an endomorphism $^{\sigma}u\ \in \End(X)$ such that
$^{\sigma}u (x)=\sigma u(\sigma^{-1}x)$ for $x \in X(K_a)$ and we
get the  group homomorphism
$$\kappa_{X}: \Gal(K) \to \Aut(\End(X));
\quad\kappa_{X}(\sigma)(u)=\ ^{\sigma}u \quad \forall \sigma \in
\Gal(K),u \in \End(X).$$ It is well-known that $\End_K(X)$
coincides with the subring of $\Gal(K)$-invariants in $\End(X)$,
i.e., $\End_K(X)=\{u\in \End(X)\mid\  ^{\sigma}u\ =u \quad \forall
\sigma \in \Gal(K)\}$. It is also well-known that $\End(X)$
(viewed as a group with respect to addition) is a free commutative
group of finite rank and $\End_K(X)$ is its {\sl pure} subgroup,
i.e., the quotient $\End(X)/\End_K(X)$ is also  a free commutative
group of finite rank.

The following observation plays a crucial role in this paper.

\begin{rem}
\label{alice}
It is known  \cite{Silverberg} that all the
endomorphisms of $X$ are defined over $K(X_{\ell^2})$; in
particular,
$$\Gal(K(X_{\ell^2})) \subset \ker(\kappa_{X})\subset \Gal(K).$$
 This implies that if
$\Gamma_K: =\kappa_{X}(\Gal(K)) \subset \Aut(\End(X))$ then there
exists a surjective homomorphism
$\kappa_{X,\ell^2,K}:\tilde{G}_{\ell^2,X,K} \twoheadrightarrow
\Gamma_K$ such that the composition
$$\Gal(K)\twoheadrightarrow \Gal(K(X_{\ell^2})/K)=
\tilde{G}_{\ell^2,X,K}\stackrel{\kappa_{X,\ell^2,K}}{\twoheadrightarrow}
\Gamma_K$$ coincides with $\kappa_X$. Clearly,
$$\End_K(X)=\End(X)^{\Gamma_K}.$$
\end{rem}

Clearly,  $\End(X)$ leaves invariant the subgroup $X_{\ell}\subset
X(K_a)$. It is well-known that $u\in \End(X)$ kills $X_{\ell}$
(i.e. $u(X_{\ell})=0$) if and only if $u \in \ell\cdot\End(X)$.
This gives us a natural embedding
$$\End_K(X)\otimes\Z/\ell\Z\subset \End(X)\otimes\Z/\ell\Z\
\hookrightarrow \End_{\F_{\ell}}(X_{\ell}).$$ Clearly, the image
of $\End_K(X)\otimes\Z/\ell\Z$ lies in the centralizer of the
Galois group, i.e., we get an embedding
$$\End_K(X)\otimes\Z/\ell\Z
\hookrightarrow
\End_{\Gal(K)}(X_{\ell})=\End_{\tilde{G}_{\ell,X,K}}(X_{\ell}).$$

 It is well-known
 \cite[\S 21]{MumfordAV},  that $\End^0(X)$ is a semisimple
finite-dimensional $\Q$-algebra.  Clearly, the natural map
$\Aut(\End(X)) \to \Aut(\End^0(X))$ is an embedding. This allows
us to view $\kappa_{X}$ as a homomorphism
$$\kappa_{X}: \Gal(K) \to \Aut(\End(X))\subset \Aut(\End^0(X)),$$
whose image coincides with $\Gamma_K\subset \Aut(\End(X))\subset
\Aut(\End^0(X))$.
 Clearly, the subalgebra $\End^0(X)^{\Gamma_K}$  of
 $\Gamma_K$-invariants coincides with $\End_K(X)\otimes\Q$.

\begin{rem}
\label{split}
\begin{itemize}
\item[(i)]
 Let us split
the semisimple $\Q$-algebra $\End^0(X)$ into a finite direct
product $\End^0(X)= \prod_{s\in \I} D_s$
 of simple $\Q$-algebras $D_s$. (Here $\I$ is identified with the
 set of minimal two-sided ideals in $\End^0(X)$.)
Let $e_s$ be the identity element of $D_s$. One may view $e_s$ as
an idempotent in $\End^0(X)$. Clearly,
$$1_X=\sum_{s\in \I} e_s\in \End^0(X), \quad e_s e_t=0 \ \forall s\ne t.$$
There exists a positive integer $N$ such that all $N \cdot e_s$
lie in $\End(X)$. We write $X_s$ for the image $X_s:=(Ne_s) (X)$;
it is an abelian subvariety in $X$ of positive dimension. Clearly,
the sum map
$$\pi_X:\prod_s X_s \to X, \quad (x_s) \mapsto \sum_s x_s$$
is an isogeny. It is also clear that the intersection $D_s\bigcap
\End(X)$ leaves $X_s \subset X$ invariant. This gives us a natural
identification $D_s \cong \End^0(X_s)$.  One may easily check that
each $X_s$ is isogenous to a self-product of  (absolutely) simple
abelian variety. Clearly, if $s\ne t$ then $\Hom(X_s,X_t)=0$.
\item[(ii)]
We write $C_s$ for the center of $D_s$. Clearly, $C_s$ coincides
with the center of $\End^0(X_s)$ and
the center $C$ of $\End^0(X)$ coincides with $\prod_{s\in
\I}C_s=\oplus_{s\in S}C_s$.
\item[(iii)]
Obviously, all the sets $$\{e_s\mid s\in \I\}\subset \oplus_{s\in
\I}\Q\cdot e_s\subset \oplus_{s\in \I}C_s=C$$ are stable under the
Galois action $\Gal(K)
\stackrel{\kappa_X}{\longrightarrow}\Aut(\End^0(X))$. In
particular, there is a continuous homomorphism from $\Gal(K)$ to
the group $\Perm(\I)$ of permutations of $\I$ such that its kernel
contains $\ker(\kappa_X)$ and
$$e_{\sigma(s)}=\kappa_X(\sigma)(e_s)=\ ^{\sigma}e_s ,
\ ^{\sigma}(C_s)\ =C_{\sigma(s)}, \ ^{\sigma}(D_s)=D_{\sigma(s)}
 \quad \forall \sigma\in \Gal(K), s\in \I.$$
 By Remark \ref{alice}, the homomorphism $\Gal(K)\to\Perm(\I)$ factors
 through $\tilde{G}_{X,\ell^2,K}$.
It follows that $X_{\sigma(s)}=Ne_{\sigma(s)}(X)=\sigma (Ne_s
(X))=\sigma(X_s)$; in particular, abelian subvarieties $X_s$ and
$X_{\sigma(s)}$ have the same dimension. Clearly, $u\mapsto\
^{\sigma}u$ gives rise to an isomorphism of $\Q$-algebras
$\End^0(X_{\sigma(s)})\cong \End^0(X_s)$.
\item[(iv)]
 Assume that $\End_K(X)$ has no zero divisors. Then the following conditions hold
 \cite[remark 1.4(iv)]{ZarhinL}:
 $\I$  consists of one Galois orbit.
 In particular, $\tilde{G}_{X,\ell^2,K}$
acts transitively on $\I$ and therefore has a subgroup of index $\#(I)$.
\end{itemize}
\end{rem}

In the next section we will need the following criterion of simplicity for algebras over fields.

\begin{lem}
\label{herst}
Let $F$ be a field, $R$ a finite-dimensional central
simple $F$-algebra. Let $A\subset R$ be a simple $F$-subalgebra with
the same identity element. (In particular, $A$ contains the center
$F$ of $R$.) Then the centralizer $C_R(A)$ of $A$ in $R$ is a simple
$F$-algebra.
\end{lem}

\begin{proof}
This is a special case of Theorem 4.3.2 on p. 104 of
\cite{Herstein}
\end{proof}

\section{Abelian varieties with multiplications}
\label{abmult} Let $E$ be a number field. Let $(X, i)$ be a pair
consisting of an abelian variety $X$ of positive dimension over
$K_a$ and an embedding $i:E \hookrightarrow  \End^0(X)$. Here
$1\in E$ must go to $1_X$.

\begin{rem}
\label{reldeg}
 It is well known  \cite[Prop. 2 on p.
36]{Shimura2}) that the degree $[E:\Q]$ divides $2\dim(X)$, i.e.
$$d_{X,E}:=\frac{2\dim(X)}{[E:\Q]}$$
is a positive integer.
\end{rem}

Let us denote by $\End^0(X,i)$ the centralizer of $i(E)$ in
$\End^0(X)$. Clearly, $i(E)$ lies in the center of the
finite-dimensional $\Q$-algebra $\End^0(X,i)$. It follows that
$\End^0(X,i)$ carries the natural structure of a finite-dimensional
$E$-algebra.

\begin{rem}
\label{ss} The $E$-algebra $\End^0(X,i)$ is semisimple. Indeed, in
the notation of Remark \ref{split}, $\End^0(X)= \prod_{s\in\I} D_s$
where all  $D_s=\End^0(X_s)$ are simple $\Q$-algebras. If
$\pr_s:\End^0(X) \twoheadrightarrow D_s$ is the corresponding
projection map and
  $D_{s,E}$ is the centralizer of $\pr_s i(E)$ in $D_s$ then one may easily check that
$\End^0(X,i)=\prod_{s\in\I} D_{s,E}$. Clearly, $\pr_s i(E)\cong E$
is a simple $\Q$-algebra.
We write $i_s$ for the composition $ \pr_s i: E \hookrightarrow
\End^0(X)\twoheadrightarrow D_{s} \cong \End^0(X_s)$. We have
$D_{s,E}=\End^0(X_s,i_s)$. Let $C_s$ be the center of $D_s$. Clearly,
$C_s$ is a field.

If $C_s$ is contained in $\pr_s i(E)\cong E$ then
 $\pr_s i(E)\cong E$ is a simple $C_s$-algebra. Now it follows from
 Lemma \ref{herst}
that $D_{s,E}$ is also a
{\sl simple} $C_s$-algebra.

Now suppose that $C_s$ is not contained in $\pr_s i(E)$. Let us
consider the natural homomorphism
$$C_s\otimes_{\Q} E \to D_s,$$
whose image $C_s E$ is a commutative semisimple $C_s$-subalgebra of
$D_s=\End^0(X_s)$. Clearly, $D_{s,E}$  coincides with the
centralizer of $C_s E$ in $D_s$. By Lemma \ref{herst}, if $C_s E$ is
a field then $D_{s,E}$ is a simple $C_s$-algebra and therefore is
also a simple $E$-algebra. In general case $C_s E$ splits into a
finite direct sum
$$C_s E= \oplus_{j\in J} F_j$$ of fields $F_j$. Let $e_{s,j}$ be the
identity element of $F_j$. Then $\sum_{j\in J}e_{s,j}=1_{X_s}$ and
$e_{s,j}e_{s,j^{\prime}}=0$ for each  $j, j^{\prime} \in J$ with $j
\ne j^{\prime}$.

There exists a positive integer $M$ such that all $M \cdot e_{s,j}$
lie in $\End(X_s)$. We write $Y_j$ for the image $Y_j:=(M e_{s,j})
(X_s)$; it is an abelian subvariety in $X_s$ of positive dimension
provided with natural embedding
$$i_{s,j}: F_j \hookrightarrow \End^0(Y_j)$$
that sends $e_{s,j}$ to the identity automorphism of $Y_j$. Let us
consider the following homomorphisms of abelian varieties
 $$\psi_s: \prod _{j\in J}Y_j \to X_s, \ \{y_j\}_{j\in J} \mapsto \sum_{j\in J}
 y_j,$$
 $$\phi_s: X_s \to \prod _{j\in J}Y_j, \ x \mapsto \{(M e_{s,j})x\}_{j\in
 J}.$$
 Clearly, $\psi_s\phi_s=M \cdot 1_{X_s}$ while $\phi_s \psi_s$ is
 multiplication by $M$ in $\prod _{j\in J}Y_j$. In particular, both
 $\psi_s$ and $\phi_s$
 are isogenies of abelian varieties over $K_a$. We also have
$$D_{s,E}=\psi_s[\oplus_{j\in J}\End^0(Y_j, i_{s,j})]\psi_s^{-1}\cong \oplus_{j\in J}\End^0(Y_j, i_{s,j}).$$
Since $X_s$ is isogenous to a self-product of an absolutely simple
abelian variety, say, $Z_s$, the abelian subvariety $Y_j$ is also
isogenous to a  self-product of $Z_s$. In particular, $\End^0(X_s)$
and $\End^0(Y_j)$ have the same center, namely, the field $C_s$.
Since $F_j$ contains $C_s$,   Lemma \ref{herst} implies that the centralizer $\End^0(Y_j, i_{s,j})$ (of  $F_j$)
is a simple $F_j$-algebra and therefore is a simple $E$-algebra.

This implies  that $D_{s,E}$ is a {\sl semisimple} $E$-algebra and
therefore $\End^0(X,i)$ is a semisimple $E$-algebra.

Recall that 
$D_{s,E}=\End^0(X_s,i_s)$ and therefore
$$\End^0(X,i)=\prod_{s\in\I}  D_{s,E}=   \prod_{s\in\I} \End^0(X_s,i_s) \eqno(1) .$$
Now it is clear that $\End^0(X,i)$ is a simple $E$-algebra if and only if
$\End^0(X)$ is a simple $\Q$-algebra (i.e., $\I$ is a singleton
$\{s\}$) and $C_s E$ is a field (i.e., the corresponding $J$ is a
singleton).
\end{rem}

Let $\OC$ be the ring of integers in $E$. Let $\a$ be a non-zero
ideal in $\OC$. It is well known that the quotient $\OC/\a$ is a
finite commutative ring.

Let $\lambda$ be a maximal ideal in $\OC$. We write $k(\lambda)$
for the corresponding (finite) residue field $\OC/\lambda$ and
$\ell$ for $\fchar(k(\lambda))$.  We have
 $$\OC\supset\lambda\supset \ell\cdot\OC.$$

\begin{rem}
\label{totram} Let us assume  that $\lambda$ is the only maximal
ideal of $\OC$ dividing $\ell$, i.e., $\ell\cdot\OC=\lambda^b$
where the positive integer $b$ satisfies $[E:\Q]=b \cdot
[k(\lambda):\F_{\ell}]$.

\begin{itemize}
\item[(i)] We have
 $\OC\otimes\Z_{\ell}=\OC_{\lambda}$ where $\OC_{\lambda}$ is
the completion of $\OC$ with respect to $\lambda$-adic topology.
It is well known that  $\OC_{\lambda}$ is a local principal ideal
domain, its only maximal ideal is $\lambda\OC_{\lambda}$ and
$$k(\lambda)=\OC/\lambda=\OC_{\lambda}/\lambda\OC_{\lambda}.$$
One may easily check that
$\ell\cdot\OC_{\lambda}=(\lambda\OC_{\lambda})^b$.

 \item[(ii)]
 Let us choose an element $c \in \lambda$ that does not
lie in $\lambda^2$. Clearly,
$$\lambda\OC_{\lambda}=c\cdot\OC_{\lambda}, \ \ell \cdot\OC_{\lambda}=c^b \cdot\OC_{\lambda}.$$
This implies that there exists a unit $u\in\OC_{\lambda}^*$ such
that $\ell=u c^b$. It follows from the unique factorization of
ideals in $\OC$ that
$$\lambda=\ell\cdot\OC + c\cdot\OC.$$
\item[(iii)] Notice that
$$E_{\lambda}=E\otimes_{\Q}\Q_{\ell}=\OC\otimes\Q_{\ell}=
\OC_{\lambda}\otimes_{\Z_{\ell}}\Q_{\ell}$$ is the field
coinciding with the completion of $E$ with respect to
$\lambda$-adic topology.
\end{itemize}
\end{rem}

Suppose that $X$ is defined over $K$ and $i(\OC) \subset
\End_K(X)$. Then we may view elements of $\OC$ as
$K$-endomorphisms of $X$. We write $\End(X,i)$ for the centralizer
of $i(\OC)$ in $\End(X)$ and $\End_K(X,i)$ for the centralizer of
$i(\OC)$ in $\End_K(X)$. Obviously,  $\End(X,i)$ is a pure
subgroup in $\End(X)$ and $\End_K(X,i)=\End(X,i)\bigcap \End_K(X)$
is a pure subgroup in $\End(X,i)$, i.e. the quotients
$\End(X)/\End(X,i)$ and
 \newline
 $\End(X,i)/\End_K(X,i)$ are torsion-free. We have
$$i(\OC)\subset \End_K(X,i)\subset\End(X,i)\subset \End(X).$$
It is also clear that
$$\End(X,i)=\End^0(X,i)\bigcap \End(X),$$
$$\End^0(X,i)=\End(X,i)\otimes\Q\subset
\End(X)\otimes\Q=\End^0(X).$$
 Clearly, $\End(X,i)$ carries the natural
structure of a finitely generated torsion-free $\OC$-module. Since
$\OC$ is a Dedekind ring, it follows that there exist non-zero
ideals $\b_1,\ldots , b_t$ in $\OC$ such that
$$\Z \cdot 1_X\subset \End(X,i)\cong \b_1\oplus \cdots \oplus \b_t$$
(as $\OC$-modules).

Clearly, $\kappa(\sigma)(\End(X,i))=\End(X,i)$ for all $\sigma\in
\Gal(K)$ and
$$\End_K(X,i)=\{u\in \End(X,i)\mid\  ^{\sigma}u\ =u \quad \forall
\sigma \in \Gal(K)\}.$$

\begin{lem}
\label{transitive} Suppose that $X$ is defined over $K$ and
$i(\OC) \subset \End_K(X)$. If $\ \End_K(X,i)$ has no zero
divisors then, in the notations of Remark \ref{split}, $\I$  consists
of one Galois orbit; in particular, all $X_s$ have the same
dimension equal $\dim(X)/\#(\I)$.
\end{lem}

\begin{proof}
Clearly, $\End_K(X,i)\otimes\Q$ also has no zero divisors.
If $\J\subset \I$ is a non-empty Galois-invariant subset then $e_{\J}:=\sum_{s\in\J}e_s$
is a non-zero element of $\End_K(X,i)\otimes\Q$. Now, if $\I$ contains two non-empty disjoint
Galois-invariant subsets $\J_1$ and $\J_2$ then $e_{\J_1} e_{\J_2}=0$.
\end{proof}

\begin{thm}
\label{multtangent} Suppose that $E$ is a number field that is
normal over $\Q$. Suppose that $K$ is a field of characteristic
zero that contains a subfield isomorphic to $E$. Suppose that $X$
is defined over $K$ and $i(\OC) \subset \End_K(X)$. Let
$\Lie_K(X)$ be the tangent space to $X$ at the origin, which
carries the natural structure of a $E\otimes_{\Q}K$-module. For each
field embedding $\tau:E\hookrightarrow K$ we put
$$\Lie_K(X)_{\tau}=\{z\in \Lie_K(X)\mid i(e)z=\tau(e)z \ \forall e\in E\},$$
$$n_{\tau}(X,i)=\dim_K(\Lie_K(X)_{\tau}).$$
If $\End_K(X,i)$ has no zero divisors then every $n_{\tau}$ is divisible by $\#(\I)$.
In particular, if the greatest common divisor of all $n_{\tau}$ is $1$ then $\#(\I)=1$,
 i.e., 
  $\End^0(X)$
  is a simple $\Q$-algebra, and $X$ is isogenous over $K_a$
 to a self-product of an absolutely simple abelian variety.
\end{thm}
\begin{rem}
The nonnegative integers $n_{\tau}(X,i)$ do not depend on the
choice of $K$.
\end{rem}

\begin{proof}[Proof of Theorem \ref{multtangent}]
The embedding $i$ gives rise to the embeddings
$$i_s:E \hookrightarrow \End^0(X_s), \ i_s(1)=1_{X_s}, i_s(\OC)\subset \End(X_s)$$
and
$$\Lie_{K_a}(X)_{\tau}=\oplus_{s\in\I}\Lie_{K_a}(X_s)_{\tau}.$$
On the other hand, if $s=\sigma(t)$ for $s,t\in\I, \sigma\in\Gal(K)$
then
$$\Lie_{K_a}(X_s)_{\tau}=\sigma (\Lie_{K_a}(X_t)_{\tau}).$$
Using Lemma \ref{transitive}, we conclude that for each $\tau$ (and $s$)
$$\dim_{K_a}(\Lie_{K_a}(X)_{\tau})=\#(\I)\cdot \dim_{K_a}(\Lie_{K_a}(X_s)_{\tau})$$
is divisible by $\#(I)$. In order to finish the proof, one has only to notice that
$$\Lie_{K_a}(X)_{\tau}=\Lie_{K}(X)_{\tau}\otimes_K K_a,$$
since $E$ is normal and $K$ contains the subfield isomorphic to $E$,
\end{proof}

\begin{lem}
\label{nocommon} In the notation and assumptions of Theorem
\ref{multtangent}, let us assume that $\End^0(X,i)$ is a central
simple $E$-algebra. Let us define the positive integer
 $m$ as the square root
of $\dim_E(\End^0(X,i))$. Then all $n_{\tau}(X,i)$ are divisible by $m$. In particular,
 if the greatest common divisor of all $n_{\tau}(X,i)$ is $1$ then $m=1$, i.e.,
$\End^0(X,i)=i(E)\cong E$.
\end{lem}

\begin{proof}
We may assume that $K=K_a$ is algebraically closed. Then for each
field embedding $\tau:E\hookrightarrow K$ the $K$-algebra
$\End^0(X,i)\otimes_{E,\tau}K$ is isomorphic to the matrix algebra
$\M_m(K)$. On the other hand, every $\Lie_{K}(X)_{\tau}$ carries the 
natural structure of a $\End^0(X,i)\otimes_{E,\tau}K$-module. It
follows that the $K$-dimension of $\Lie_{K}(X)_{\tau}$
 is divisible by $m$.
\end{proof}

 Let us assume
that $\fchar(K)$ does not divide the order of  $\OC/\a$ and put
$$X_{\a}:=\{x \in X(K_a)\mid i(e)x=0 \quad \forall e\in
\a\}.$$
 For example,  assume that $\ell \ne \fchar(K)$.  Then the order
 of $\OC/\lambda=k(\lambda)$ is a power of $\ell$ and
$$X_{\lambda}\subset X_{\ell}.$$ Assume in addition that $\lambda$ is the only
maximal ideal of $\OC$ dividing $\ell$ and pick
$$c \in \lambda \setminus \lambda^2\subset \lambda \subset
 \OC.$$ By Remark \ref{totram}(ii),  $\lambda$ is
generated by $\ell$ and $c$. Therefore
$$X_{\lambda}=\{x\in X_{\ell}\mid cx=0\}\subset X_{\ell} \ \eqno(2).$$

 Clearly, $X_{\a}$ is a Galois submodule of $X(K_a)$. It is
also clear that $X_{\a}$ carries the natural structure of a
$\OC/\a$-module. (It is known \cite[Prop. 7.20]{Shimura} that this
module is free.)

Obviously, every endomorphism from $\End(X,i)$ leaves invariant
the subgroup $X_{\a} \subset X(K_a)$ and induces an endomorphism
of the $\OC/\a$-module $X_{\a}$. This gives rise to a natural
homomorphism
$$\End(X,i) \to \End_{\OC/\a}(X_{\a}),$$
whose kernel contains $\a\cdot\End(X,i)$. We claim that actually,
the kernel coincides with $\a\cdot\End(X,i)$, i. e. there is an
embedding
$$\End(X,i)\otimes_{\OC}\OC/\a \hookrightarrow \End_{\OC/\a}X_{\a}\ \eqno(3).$$
In order to prove it, let us denote by $n$ the order of $\OC/\a$.
 We have $n\OC\subset\a$.  Since $\OC$ is a Dedekind ring, there
exists a non-zero ideal $\a'$ in $\OC$ such that $\a'\a=n\OC$.
Clearly
$$\a'X_n\subset X_{\a}.$$
It follows that if $u\in\End(X,i)$ kills $X_{\a}$ then $hu=uh$
kills $X_n$ for each $h\in\a'$. This implies that
$$\a' u\subset n\End(X)\bigcap\End(X,i)=n\End(X,i)= n(\b_1\oplus \cdots \oplus \b_t)=
n\b_1\oplus \cdots \oplus n\b_t.$$ Since
$$n\b_1\oplus \cdots \oplus n\b_t=\a'\a\b_1\oplus \cdots \oplus \a'\a\b_t,$$ we have
 $u \in\a\b_1\oplus \cdots \oplus \a\b_t= \a\cdot\End(X,i)$ and we are
 done.

Now we concentrate on the case of $\a=\lambda$, assuming that
$$\ell \ne \fchar(K).$$
Then $X_{\lambda}$ carries the natural
structure of a $k(\lambda)$-vector space provided with the
structure of Galois module and (3) gives us the embedding
$$\End(X,i)\otimes_{\OC}k(\lambda)\hookrightarrow\End_{k(\lambda)}(X_{\lambda})
. \eqno(4)$$ Further we will identify
$\End(X,i)\otimes_{\OC}k(\lambda)$ with its image in
$\End_{k(\lambda)}(X_{\lambda})$. We write
$$\tilde{\rho}_{\lambda,X}:\Gal(K) \to
\Aut_{k(\lambda)}(X_{\lambda})$$ for the corresponding
(continuous) homomorphism defining the Galois action on
$X_{\lambda}$. It is known \cite[Prop. 7.20] {Shimura} (see also
\cite{Ribet2}) that
$$\dim_{k(\lambda)}X_{\lambda}= \frac{2\dim(X)}{[E:\Q]}:=d_{X,E}.$$

 Let us put
$$\tilde{G}_{\lambda,X}=\tilde{G}_{\lambda,X,K}:=\tilde{\rho}_{\lambda,X}(\Gal(K)) \subset
\Aut_{k(\lambda)}(X_{\lambda}).$$ Clearly, $\tilde{G}_{\lambda,X}$
coincides with the Galois group of the field extension
$K(X_{\lambda})/K$ where $K(X_{\lambda})$ is the field of
definition of all points in $X_{\lambda}$.

It is also clear that the image of $\End_K(X,i)
\otimes_{\OC}k(\lambda)$ lies in the centralizer
$\End_{\tilde{G}_{\lambda,X,K}}(X_{\lambda})$ of $\tilde{G}_{\lambda,X,K}$ in
$\End_{k(\lambda)}(X_{\lambda})$. This implies the following
statement.

\begin{lem}
\label{overK} Suppose that $i(\OC)\subset \End_K(X)$. If $\lambda$
is a maximal ideal in $\OC$ such that $\ell \ne \fchar(K)$ and
  $\End_{\tilde{G}_{\lambda,X,K}}(X_{\lambda})=k(\lambda)$
then $\End_K(X,i)=\OC$.
\end{lem}
\begin{proof}
Clearly, $\End_K(X,i)$ is a finitely generated torsion free
$\OC$-module and therefore is isomorphic to a direct sum of say,
$r$ (non-zero) ideals in $\OC$. Clearly, $r=1$ and therefore
$\OC\subset\End_K(X,i)\subset E$. Since $\OC$ is the maximal order
in $E$, we conclude that $\End_K(X,i)=\OC$.
\end{proof}

Since $X_{\lambda}\subset X_{\ell}\subset X_{\ell^2}$, we get the inclusion
$K(X_{\lambda})\subset K(X_{\ell^2})$, which gives rise to the
surjective group homomorphism
$$\tau_{\ell^2,\lambda}:\tilde{G}_{\ell^2,X,K}\twoheadrightarrow
\tilde{G}_{\lambda,X}.$$

\begin{rem}
\label{minimal}
Recall \cite{FT} (see also \cite[p. 199]{ZarhinP})
that a surjective homomorphism of finite groups
$\pi:\GG_1\twoheadrightarrow \GG$ is called a {\sl minimal cover}
if no proper subgroup of $\GG_1$ maps onto $\GG$. One may easily check
that if $H$ is a subgroup in $\GG_1$ then $(\GG:\pi(H))$ divides $(\GG_1:H)$
(and the ratio is the index of $\ker(\pi)\bigcap H$ in $\ker(\pi)$).
Also, if $H$ is normal in $\GG_1$ then $\pi(H)$ is normal in $\GG$.
\end{rem}

The following notion is a natural generalization of Definition 2.2 in \cite{ElkinZarhin}.

\begin{defn}
Suppose that $\lambda$ is the only maximal ideal of $\OC$ dividing
$\ell$. We say that $K$ is $\lambda$-balanced with respect to $X$
if $\ell \ne \fchar(K)$ and
$$\tau_{\ell^2,\lambda}:
\tilde{G}_{\ell^2,X,K}\twoheadrightarrow \tilde{G}_{\lambda,X,K}$$
 is a minimal cover.
\end{defn}

\begin{rem}
\label{bigenough}
 There always exists a subgroup $H
\subset \tilde{G}_{\ell^2,X,K}$ such that $H\to\tilde{G}_{\lambda,X,K}$ is
surjective and a minimal cover. (Indeed, one has only to choose a subgroup
$H\subset\tilde{G}_{\ell^2,X,K}$ of smallest possible order
such that $\tau_{\ell^2,\lambda}(H)=\tilde{G}_{\lambda,X,K}$.)
Let us put $L=K(X_{\ell^2})^H$. Clearly,
$$K \subset L \subset K(X_{\ell^2}), \ L\bigcap K(X_{\lambda})=K$$
and $L$ is a maximal overfield of $K$ that enjoys these
properties. It is also clear that
$$K(X_{\lambda})\subset L(X_{\lambda}),\ L(X_{\ell^2})=K(X_{\ell^2}), \ H=\tilde{G}_{\ell^2,X,L},
\tilde{G}_{\lambda,X,L}=\tilde{G}_{\lambda,X,K}$$
 and $L$ is $\lambda$-{\sl balanced}  with respect to $X$.
\end{rem}

In order to describe $\tilde{\rho}_{\lambda,X}$ explicitly, let us
assume for the sake of simplicity that $\lambda$ is the only
maximal ideal of $\OC$ dividing $\ell$. We write $b$ for the
corresponding ramification index and pick $c\in \lambda\setminus
\lambda^2$.

Let $T_{\ell}(X)$ be the $\Z_{\ell}$-Tate module of $X$. Recall
that $T_{\ell}(X)$ is a free $\Z_{\ell}$-module of rank $2\dim(X)$
provided with the continuous action
$$\rho_{\ell,X}:\Gal(K) \to \Aut_{\Z_{\ell}}(T_{\ell}(X))$$
and the natural embedding
$$\End_K(X)\otimes\Z_{\ell} \hookrightarrow \End_{\Z_{\ell}}(T_{\ell}(X)),$$
whose image commutes with $\rho_{\ell,X}(\Gal(K))$. In particular,
$T_{\ell}(X)$ carries the natural structure of an
$\OC\otimes\Z_{\ell}=\OC_{\lambda}$-module; it is known
\cite[Prop. 2.2.1 on p. 769]{Ribet2} that the
$\OC_{\lambda}$-module $T_{\ell}(X)$ is free of rank $d_{X,E}$.
There is also the natural isomorphism of Galois modules
$$X_{\ell}=T_{\ell}(X)/\ell T_{\ell}(X),$$
which is also an isomorphism of $\End_K(X)\supset \OC$-modules.
Applying (2) and Remark \ref{totram}(ii), we conclude  that
 the $\OC[\Gal(K)]$-module $X_{\lambda}$ coincides with
$$c^{-1}\ell T_{\ell}(X)/\ell T_{\ell}(X)=c^{b-1}T_{\ell}(X)/c^b
T_{\ell}(X)=T_{\ell}(X)/c T_{\ell}(X)= $$ $$T_{\ell}(X)/c
\OC_{\lambda} T_{\ell}(X) =
T_{\ell}(X)/(\lambda\OC_{\lambda})T_{\ell}(X).$$ Hence
$$X_{\lambda}=T_{\ell}(X)/(\lambda\OC_{\lambda})T_{\ell}(X)=T_{\ell}(X)\otimes_{\OC_{\lambda}}k(\lambda).$$

\begin{thm}
\label{mainAV} Suppose that $E$ is a number field that is normal
over $\Q$. Suppose that $K$ is a field of characteristic zero that
contains a subfield isomorphic to $E$. Suppose that $X$ is defined
over $K$ and $i(\OC) \subset \End_K(X)$. We write $n_{X,E}$ for
the greatest common divisor of all $n_{\tau}(X,i)$ where $\tau$
runs through the set of embeddings $\tau:E\hookrightarrow K$.
Suppose that $\lambda$ is the only maximal ideal in $\OC$ that
lies above prime $\ell$, the group $\tilde{G}_{\lambda,X}$ does
not contain a  subgroup, whose index divides $n_{X,E}$, except
$\tilde{G}_{\lambda,X}$ itself, and
$\End_{\tilde{G}_{\lambda,X}}(X_{\lambda})=k(\lambda)$. Then:
\begin{itemize}
\item[(i)]
 $\End^0(X)$ is a simple $\Q$-algebra, i.e., $X$ is isogenous over $K_a$
 to a self-product of an absolutely simple abelian variety.
\item[(ii)]
Suppose that $\tilde{G}_{\lambda,X}$ does not contain a  normal subgroup, whose index divides
$d_{X,E}$ except
$\tilde{G}_{\lambda,X}$ itself. Assume also that  $n_{X,E}=1$.
Then:

$\End^0(X,i)=i(E)\cong E$ and $\End(X,i)=i(\OC)\cong\OC$.
\end{itemize}
\end{thm}

\begin{proof}
Using Remark \ref{bigenough}, we may assume (extending if necessary the ground field)
that $K$ is $\lambda$-balanced with respect to $X$. In particular,
there is no proper subgroup
$H\subset \tilde{G}_{\ell^2,X,K}$ with
 $\tau_{\ell^2,\lambda}(H)=\tilde{G}_{\lambda,X,K}$.
Since $\End_{\tilde{G}_{\lambda,X}}(X_{\lambda})=k(\lambda)$, it
follows from Lemma \ref{overK} that $\End_K(X,i)=i(\OC)\cong\OC$
and therefore has no zero divisors. It follows from Lemma
\ref{transitive} that $\I$  consists of one Galois orbit. Assume
that $\I$ is not a singleton. By Remark \ref{split},
$\tilde{G}_{\ell^2,X,K}$ contains a subgroup $H$ of index
$\#(\I)$. It follows from Remarks \ref{bigenough} and
\ref{minimal} that
$H_0:=\tau_{\ell^2,\lambda}(H)\ne\tilde{G}_{\lambda,X,K}$; in
addition the index of $H_0$ in $\tilde{G}_{\lambda,X,K}$ divides
$\#(\I)$. By Theorem \ref{multtangent}, $\#(\I)$ divides
$n_{X,E}$. It follows that  the index of $H_0$ in
$\tilde{G}_{\lambda,X,K}$ divides $n_{X,E}$. But such $H_0$ does
not exist: the obtained contradiction proves that $\I$ is  a
singleton, i.e., $X=X_s$, $\End^0(X)=\End^0(X_s)$ is a simple
$\Q$-algebra. 
This proves (i).

Let us prove (ii).
Assume for a while that $\End^0(X,i)$ is a simple
$E$-algebra. Then the center $C$ of the simple $E$-algebra
$\End^0(X,i)$ is a field that contains $i(E)\cong E$. By Remark
\ref{reldeg},  $[C:\Q]$ divides $2\dim(X)$. This implies that
$[C:E]$ divides $d_{X,E}$. Suppose that $C\ne E$, i.e., $[C:E]>1$.
Since $\End_K(X,i)=\OC$, we have
$$E=\End^0_K(X,i)=C^{\Gal(K)}=E^{\tilde{G}_{\ell^2,X,K}},$$
i.e., $C/E$ is a Galois extension and  there is a surjective
homomorphism $\tilde{G}_{\ell^2,X,K}\twoheadrightarrow \Gal(C/E)$,
whose kernel $H_1$ is a normal subgroup in $\tilde{G}_{\ell^2,X,K}$
of index $[C:E]>1$. Clearly, $H_2:=\tau_{\ell^2,\lambda}(H_1)$ is a
normal subgroup in $\tilde{G}_{\lambda,X,K}$ and its index is
greater than $1$ but divides $[C:E]$. Since $[C:E]$ divides
$d_{X,E}$, the index of $H_2$ in $\tilde{G}_{\lambda,X,K}$ divides
$d_{X,E}$. But such $H_2$ does not exist: the obtained contradiction
proves that $[C:E]=1$, i.e., $C=E$ and $\End^0(X,i)$ is a {\sl
central} simple $E$-algebra.  It follows from Lemma \ref{nocommon}
that $\End^0(X,i)=i(E)\cong E$. Since $\OC\cong i(\OC)\subset
\End(X,i)\subset i(E)\cong E$ and $\OC$
 is the maximal order in $E$, we conclude that $\OC\cong i(O)= \End(X,i)$.

 So, we are done, under an additional assumption that $\End^0(X,i)$ is a simple
$E$-algebra. Now, let us do the general case. Further we will
identify $E$ with its image in $\End^0(X,i)\subset \End^0(X)$. We
are going to prove that, in the notation of Remark \ref{ss}, the set
$J$ is a singleton: this would imply that $\End^0(X,i)$ is a simple
$E$-algebra. By (i), $X=X_s$. Let $C_s$ be the center of
$\End^0(X_s)$. Since $X=X_s$ is defined over $K$, the center $C_s$
is stable under the natural action of $\Gal(K)$. Since $E$ consists
of $K$-endomorphisms of $X$, the compositum $C_s E\subset
\End^0(X_s)$ is also stable under the action of $\Gal(K)$. Clearly,
$E$ coincides with the subalgebra of Galois invariants of $C_s E$,
because $C_s E$ commutes with $E$ and $\End_K(X,i)=\OC$.

By Remark \ref{ss}, $C_s E$ is a direct sum $\oplus_{j\in J} F_j$ of
fields $F_j$. The set of the identity elements $e_{s,j}$'s of
$F_j$'s is the set of minimal idempotents of $C_s E$ and therefore
is Galois-stable. I claim that  the Galois action on the set of
$e_{s,j}$'s is transitive. Indeed, suppose that $J_1$ and $J_2$ are
two disjoint Galois orbits in $J$. Then
$$u_1=\sum_{j\in J_1} e_{s,j}, \ u_2=\sum_{j\in J_2} e_{s,j}$$
are nonzero Galois invariants, whose product is zero. Therefore $u_1$
and  $u_2$ are nonzero elements of the field $E$, whose product is
zero, which is nonsense. The obtained contradiction proves that the
Galois action on the set of $e_{s,j}$'s is transitive. Let us consider the tangent space
$\Lie_K(X)$ to $X$ at the origin; it is a finite-dimensional vector
space over $K$. If we consider the tangent space to $X_s=X$ over
$K_a$ then we get the $K_a$-vector space
$$\Lie_{K_a}(X_s)=\Lie_K(X)\otimes_K K_a.$$
The endomorphism algebra $\End^0(X)=\End^0(X_s)$ acts by
functoriality by $K_a$-linear transformations on the tangent space
$\Lie_{K_a}(X_s)$. (Here we use that $\fchar(K_a)=\fchar(K)=0$.)

 The Galois group $\Gal(K)$ acts naturally by semi-linear
automorphisms on $\Lie_{K_a}(X_s)$. Let us consider the $K_a$-vector
subspaces
$$W_j=(Me_{s,j})\Lie_{K_a}(X_s)= e_{s,j}\Lie_{K_a}(X_s)\subset
\Lie_{K_a}(X_s).$$  Clearly, $W_j$'s are $E$-invariant,
{\bf Galois-conjugate} and
$$\Lie_{K_a}(X_s)=\oplus_{j\in J} W_j.$$
Since the action of $E$ on $\Lie_{K_a}(X_s)=\Lie_K(X)\otimes_K K_a$,
is defined over $K$, for each $a \in E$ the trace of the
corresponding $K_a$-linear operator $W_j \to W_j$ induced by $a$ does not depend
on $j$. This implies that for each field embedding $\tau: E \hookrightarrow K_a$
the corresponding multiplicity $n_{\tau} (X,i)$ is divisible by
$\#(J)$. (Recall that $\tau(E)$ always sits in $K$.) It follows that
$n(X,E)$ is divisible by $\#(J)$. Since $n(X,E)=1$, the set $J$ is a
singleton and we are done.
\end{proof}

\section{superelliptic jacobians}
\label{proofm} We keep the notation of Section \ref{one}. In
particular, $K$ is a field of characteristic zero that contains
$\zeta$: a primitive $p$th root of unity; $f(x)\in K[x]$ is a
polynomial of degree $n\ge 3$  and without multiple roots,
$\RR_f\subset K_a$ is the ($n$-element)  set of roots of $f$ and
$K(\RR_f)\subset K_a$ is the splitting field of $f$. We assume
that $p$ does not divide $n$ and write $\Gal(f)=\Gal(f/K)$ for the
Galois group $\Gal(K(\RR_f)/K)$ of $f$; it permutes roots of $f$
and may be viewed as a certain permutation group of $\RR_f$, i.e.,
as a subgroup of the group $\Perm(\RR_f)\cong\Sn$ of permutation
of $\RR_f$.

 We write $C_{f,p}$
for the superelliptic $K$-curve $y^p=f(x)$ and $J(C_{f,p})$ for
its jacobian. Recall that $J(C_{f,p})$ is an abelian variety that is
defined over $K$ and
$$\dim(J(C_{f,p}))=\frac{(n-1)(p-1)}{2}.$$
Let us denote by $d(n,p)$ the greatest common divisor of all integers $[ni/p]$
where $i$ runs through the set of integers
$$\{i\mid 0<i<p\}.$$

\begin{rem}
\label{NOD}
\begin{itemize}
\item[(i)] Suppose that $p$ divides $n-1$, i.e., $n=pk+1$ with
positive integer $k$. Then $[ni/p]=ik$ and therefore
$d(n,p)=k=(n-1)/p$. \item[(ii)] Suppose that $p$ does {\sl not}
divide $n-1$, i.e., $n=pk+r$ with nonnegative integer $k$ and
$1<r<p$. Then $p>2$, since $p$ does not divide $n$. Let us put
$j:=[p/r]\le [p/2]<p-1$. We have $[n\cdot 1/p]=k,
[n(j+1)/p]=(j+1)k+1$ and therefore $d(n,p)=1$.
\end{itemize}
\end{rem}

Let us put
$$E=\Q(\zeta_p), \ \OC=\Z[\zeta_p],\ \ell=p, \ \lambda =(1-\zeta_p)\Z[\zeta_p], \  k(\lambda)=\F_p$$
and
$$X=J(C_{f,p}),\ i:\Z[\zeta_p]\hookrightarrow \End_K(J(C_{f,p})), \ i(\zeta_p)=\delta_p.$$
We have $d_{J(C_{f,p}),\Q(\zeta_p)}=n-1$. It is known
\cite{Poonen,SPoonen} that $\Gal(f)$ is canonically isomorphic to
$\tilde{G}_{J(C_{f,p}),\lambda,K}$. It is also known
\cite{ZarhinCrelle,ZarhinM}, \cite[Theorem 4.14]{ZarhinPisa}
 that the centralizer of $\tilde{G}_{J(C_{f,p}),\lambda,K}$ in $\End_{\F_p}(J(C_{f,p})_{\lambda})$
  coincides with $\F_p$ if $\Gal(f)$ is {\sl doubly transitive}. Also, if $p>2$ and
  $\tau_i:E=\Q(\zeta_p)\hookrightarrow K$
  is an embedding that sends $\zeta_p$ to $\zeta^{-i}$ then
  $n_{\tau_i}(J(C_{f,p}),\Q(\zeta_p))=[ni/p]$ if $1\le i<p$ \cite[Remark 3.7]{ZarhinCamb}. This implies
  that $n_{J(C_{f,p}),\Q(\zeta_p)}=d(n,p)$.

\begin{proof}[Proof of Theorem  \ref{main}]
Extending the ground field, we may and will assume that $\Gal(f)=\G$.
Recall that $p>2$. Now the assertion (i) follows from Theorem \ref{mainAV}(i)
   combined with Remarks \ref{NOD} and \ref{proper}.
 It follows from Theorem \ref{mainAV}(ii)(i) combined with  Remark \ref{proper}
that $\Q[\delta_p]$ coincides with its own centralizer in  $\End^0(J(C_{f,p}))$; in particular,
$\Q[\delta_p]$ contains the center of $\End^0(J(C_{f,p}))$. It follows from
  \cite[Theorem 3.6]{ZarhinCamb} and \cite[Theorem 3.2]{ZarhinPisa} that
$\Q[\delta_p]$ coincides with the center of $\End^0(J(C_{f,p}))$.
This means that the centralizer of $\Q[\delta_p]$ in
$\End^0(J(C_{f,p}))$ is the whole $\End^0(J(C_{f,p}))$. Now
Theorem \ref{mainAV}(ii)(ii) combined with Remark \ref{proper}
implies that $\End^0(J(C_{f,p}))=\Q[\delta_p]\cong\Q(\zeta_p)$.
Since $\Z[\zeta_p]\cong \Z[\delta_p]$ is the maximal order in
$\Q(\zeta_p)$ and $\Z[\delta_p]\subset \End(J(C_{f,p}))$, we
conclude that $\Z[\delta_p]= \End(J(C_{f,p}))$.
\end{proof}

\begin{thm}
\label{mainASp} Suppose that $p$ is odd and there exists a prime $\ell\ne p$ and
positive integers $r$ and $d$ such that  $n=\ell^{2rd}$. Let $\F$
be the finite field $\F_{\ell^r}$. Let $\F^{2d}$ be the
$2d$-dimensional coordinate $\F$-vector space and $\SPA(2d,\F)$ be
the group of all symplectic affine transformation
$$v \mapsto A(v)+t, \ t\in \F^d, A\in \Sp(2d,\F)$$
of $\F^{2d}$. Suppose that one may identify $\RR_f$ and $\F^{2d}$
in such a way that $\Gal(f)$  becomes a permutation group of
$\F^{2d}$ that contains $\SPA(2d,\F)$.

Suppose that $\fchar(K)=0$  and $K$ contains a primitive $p$th
root of unity. If either $n=p+1$ or $p$ does not divide $n-1$ then
$$\End^0(J(C_{f,p}))=\Q[\delta_p]\cong\Q(\zeta_p), \ \End(J(C_{f,p}))=\Z[\delta_p]\cong\Z[\zeta_p].$$
\end{thm}

\begin{proof}
The permutation group $\G=\SPA(2d,\F)$ is doubly transitive,
because the symplectic group $\Sp(2d,\F)$ is transitive on the set
of all non-zero vectors in $\F^{2d}$. The group $T\cong \F^{2d}$
of all translations of $\F^{2d}$ is normal in $\SPA(2d,\F)$ and
$\G/T=\Sp(2d,\F)$.  Suppose that $\G$ contains a proper normal
subgroup say, $H$, whose index divides $n-1=\ell^{2rd}-1$. It
follows that $H$ contains all elements, whose order is a power of
$\ell$. This implies that $H$ contains $T$ and therefore $H/T$ is
a normal subgroup in $\Sp(2d,\F)$ that contains all elements of
$\Sp(2d,\F)$, whose order is a power of $\ell$. Since $\Sp(2d,\F)$
is generated by {\sl symplectic transvections} \cite[Theorem
3.25]{Artin} and every (non-identity)  symplectic transvection
over $\F$ has multiplicative order $\ell$, we conclude that
 $H/T=\Sp(2d,\F)$ and therefore $H=\SPA(2d,\F)$. Since $H$ is
proper, we get a contradiction. So such $H$ does not exist. Now
Theorem \ref{mainASp} follows from Theorem \ref{main}.
\end{proof}

In the case of hyperelliptic jacobians (i.e., when $p=2$), we still can get some information about
 the structure of the endomorphism algebra.

\begin{thm}
 \label{main2}
Suppose that $K$ has characteristic zero, $n$ is an odd integer that is greater or equal
 than $5$ and $\Gal(f)$ contains a doubly transitive
subgroup $\G$ that does not contain a proper subgroup, whose index divides $(n-1)/2$.
Then:
 \begin{itemize}
\item[(i)]
$\End^0(J(C_{f,2}))$ is a simple $\Q$-algebra, i.e., $J(C_{f,2})$ is isogenous over $K_a$
to a self-product of an absolutely simple abelian variety.
\item[(ii)]
If $\G$ does not contain a subgroup of index $2$ then
 $\End^0(J(C_{f,p}))$ is a central simple $\Q$-algebra.
\end{itemize}
\end{thm}

\begin{proof} Since $p=2$, we have
  $E=\Q(\zeta_2)=\Q, \ \dim(J(C_{f,2}))=(n-1)/2$ and  the assertion follows readily
  from Theorem 1.6 in \cite{ZarhinL}.
 \end{proof}

\section{Complements}
\label{comq}
Let $r>1$ be a positive integer and $q=p^r$. Let $\zeta_q\in\C$ be
a primitive $q$th root of unity, $E:=\Q(\zeta_q)\subset\C$ the
$q$th cyclotomic field and $\Z[\zeta_q]$ the ring of integers in
$\Q(\zeta_q)$. Suppose that $K$ is a field of characteristic zero
than contains a primitive $q$th root of unity $\zeta$. Let
$f(x)\in K[x]$ be a polynomial of degree $n\ge 4$ without multiple
roots. Let us assume that $p$ does {\sl not} divide $n$. We write
$C_{f,q}$ for the smooth projective model of smooth affine
$K$-curve $y^q=f(x)$. The map $(x,y)\mapsto (x,y^q)$ gives rise to
the surjective regular map $\pi:C_{f,q}\to C_{f,q/p}$. By Albanese
functoriality, $\pi$ induces the surjective homomorphism of
corresponding jacobians $J(C_{f,q})\to J(C_{f,q/p})$,  which we
still denote by $\pi$. We write $J^{(f,q)}$ for the identity
component of the kernel of
$$\pi:J(C_{f,q})\to J(C_{f,q/p}).$$
It is known \cite[pp. 355--356]{ZarhinM}, \cite[Theorem 5.13 and
Remark 5.14 on p. 383]{ZarhinPisa} that:

\begin{itemize}
\item[(i)]
 The automorphism
$$\delta_q:C_{f,q}\to C_{f,q}, \ (x,y)\mapsto (x,\zeta y)$$
gives rise to an embedding
$$\Z[\zeta_q]\hookrightarrow \End_K(J^{(f,q)})\subset \End(J^{(f,q)})$$
that sends $1$ to the identity automorphism of $J^{(f,q)}$.

\item[(ii)]
 The Galois modules
$J^{(f,q)}_{(1-\zeta_q)}$ and $J(C_{f,p})_{(1-\zeta_p)}$ are
isomorphic. Here ${(1-\zeta_q)}$ (resp. ${(1-\zeta_p)}$) is the
corresponding principal maximal ideal in $\Z[\zeta_q]$ (resp. in
$\Z[\zeta_p]$). In particular, if
$\lambda=(1-\zeta_q)\Z[\zeta_q]\subset \Z[\zeta_q]$ then
$\tilde{G}_{\lambda,J^{(f,q)},K}=\Gal(f)$; if $\Gal(f)$ is doubly
transitive then its centralizer in
$\End_{\F_p}(J^{(f,q)}_{\lambda})$ coincides with $\F_p$.

\item[(iii)] Let us consider the induced action of $E=\Q(\zeta_q)$
on $\Lie(J^{(f,q)})$. Let $i$ be a positive integer such that
$1\le i<q$ and $p$ does not divide $i$. If
$\tau_i:E=\Q(\zeta_q)\hookrightarrow K$
  is an embedding that sends $\zeta_q$ to $\zeta^{-i}$ then
  $n_{\tau_i}(J(C^{(f,q)}),E)=[ni/q]$.
  \end{itemize}

  \begin{lem}
  \label{dnq}
  Let $p$ be a prime, $r>1$ an integer and $q=p^r$. Let $n\ge 4$
  be an integer that is not divisible by $p$. Let $n=kq+c$ where
  $k$ and $c$ are nonnegative integers and $0<c<q$.

  Let  $d(n,q)$ be the greatest common divisor of all integers
$[ni/q]$ where $i$ runs through the set of positive integers $i$
that are not divisible by $p$ and strictly less than $q$. Then:

\begin{enumerate}
\item $d(n,q)$ divides $k$. In particular, if $q<n<2q$ then
$d(n,q)=1$. \item If $n-1$ is divisibly by $q$, i.e., $c=1$, then
$d(n,q)=k=(n-1)/q$.
 \item If $p$ is odd and $n-1$ is not divisible by $q$ then
$d(n,q)=1$. \item Suppose that $p=2$ and $q$ does not divide
$n-1$. Then:
\begin{itemize}
\item[(a)] $d(n,q)=1$ or $2$.
\item[(b)] If either $k$ is odd or
$c<q/2$ then $d(n,q)=1$. In particular, if $n<q/2$ then
$d(n,q)=1$.
\end{itemize}
\end{enumerate}
  \end{lem}

We prove Lemma \ref{dnq} at the end of this Section.

Combining  observations (i)-(iii) of this Section with  Theorem
\ref{mainAV},
   Lemma \ref{dnq} and Remark \ref{proper}, we obtain the following
   ``generalization"
   of Theorem \ref{main}.

\begin{thm}
 \label{mainq}
Suppose that $K$ has characteristic zero, $n\ge 4$ and $p$ is a
prime that does not divide $n$. Let $r>1$ be a positive integer,
$q=p^r$ and suppose that $K$ contains a primitive $q$th root of
unity. Suppose that
  $\Gal(f)$ contains a doubly transitive
subgroup $\G$ that enjoys the following property: if $q$ divides
$n-1$ and $n\ne q+1$ then $\G$
 does not contain a proper subgroup, whose index divides $(n-1)/q$.
If $p=2$ then we assume additionally that $\G$ does not contain a
subgroup of index $2$.
 Then:

 \begin{itemize}
\item[(i)] $\End^0(J^{(f,q)})$ is a simple $\Q$-algebra, i.e.,
$J^{(f,q)}$ is isogenous over $K_a$ to a self-product of an
absolutely simple abelian variety.

\item[(ii)] Let us assume  that  $\G$ does not contain a proper
normal subgroup, whose  index divides $n-1$. Assume also that either $n=q+1$ or  $q$ does not divide $n-1$.
If $p=2$ then we assume additionally that $n=kq+c$ with
nonnegative integers $k$ and $c<q$ such that either $k$ is odd or
$c<q/2$.

Then
 $$\End^0(J^{(f,q)})=\Q(\zeta_q), \ \End(J^{(f,q)})=\Z[\zeta_q].$$
 \end{itemize}
\end{thm}

\begin{cor}
 \label{main4q}
Suppose that $K$ has characteristic zero, $f(x)\in K[x]$ is an irreducible
 quartic polynomial, whose Galois group  is the full symmetric group $\mathbf{S}_4$.
 Let $p$ be an odd prime, $r>1$  a
positive integer and $q=p^r$.  Suppose that $K$ contains a primitive $q$th root of unity.
Then
$$\End^0(J^{(f,q)})=\Q(\zeta_q), \  \End(J^{(f,q)})=\Z[\zeta_q].$$
\end{cor}

\begin{proof}[Proof of Corollary \ref{main4q}]
The same proof as for Corollary \ref{main4} works in this case,
provided one replaces the reference to Theorem \ref{main} by a
reference to Theorem \ref{mainq}.
\end{proof}

\begin{cor}
\label{simpleq} Suppose that $\fchar(K)=0$ and $n\ge 5$.  Suppose
that $p$ is a prime that does not divide $n$. Let $r>1$ be a
positive integer and $q=p^r$. Assume also that either $n=q+1$ or
$q$ does not divide $n-1$. If $p=2$ then we assume additionally
that $n=kq+c$ with nonnegative integers $k$ and $c<q$ such that
either $k$ is odd or $c<q/2$.
Suppose that
  $\Gal(f)$ contains a doubly transitive simple non-abelian subgroup  $\G$.
Then
$$\End^0(J^{(f,q)})=\Q(\zeta_q), \  \End(J^{(f,q)})=\Z[\zeta_q].$$
\end{cor}

\begin{proof}[Proof of Corollary \ref{simpleq}]
The same proof as for Corollary \ref{simple} works in this case,
provided one replaces the reference to Theorem \ref{main} by a
reference to Theorem \ref{mainq} and notices that simple
non-abelian $\G$ does not contain a subgroup of index $2$.
\end{proof}

\begin{proof}[Proof of Lemma \ref{dnq}]
Clearly, $p$ does not divide $c$.
  We have
$[n\cdot 1/q]=k$ and therefore $d(n,q)$ divides $k$. This proves
(1). If $c=1$ then $k=(n-1)/q$ and $[ni/q]=ik$. This implies that
$d(n,q)=k=(n-1)/q$. This proves (2).

Now assume that $c>1$. Let $m$ be the smallest integer such that
$mc>q$. Clearly,
$$q/2>m\ge 2, \ (m-1)c< q<mc<2q$$ and therefore
$$(km+1)q<nm<(km+2)q.$$
It follows that $[nm/q]=km+1$ and therefore if $m$ is not
divisible by $p$ then $d(n,q)=1$.

Now let us assume that $m$ is divisible by $p$. Then $m+1$ is {\sl
not} divisible by $p$ and still $m+1<q$. Clearly,
$$c\ne q/2.$$
If $c<q/2$ then $mc=(m-1)c+c<q+q/2$ and
$$(m+1)c=mc+c<(q+q/2)+q/2=2q.$$ This implies that
$$(km+1)q<n(m+1)<(km+2)q$$
and therefore $[n(m+1)/q]=km+1$, which, in turn, implies that
$d(n,q)=1$.

If $c>q/2$ then $m=2$. If $p$ is odd, $2$ is not divisible by $p$
and (as we have already seen) $[n\cdot 2/q]=2k+1$ and therefore
$d(n,q)=1$. This proves (3).

 So, the only remaining case is when $p=2$ and $c>q/2$. Since $q\ge
 4$, we may take
 $i=3<q$. We have $n\cdot 3=3kq+3c$ and either
 $q<3c<2q$ or $2q<3c<3q$. In the former case $[n\cdot 3/q]=3k+1$
 and $d(n,q)=1$.

 In the latter case $[n\cdot 3/q]=3k+2$ and $d(n,q)$ divides both $k$ and
 $3k+2$. This implies that $d(n,q)=1$ or $2$; in addition, if $k$ is odd
then $d(n,q)=1$. This proves (4).
\end{proof}

\section{Corrigendum to \cite{ZarhinM} and \cite{ZarhinMZ}}

\cite{ZarhinM}, Remark 4.1 on p. 351: the condition $n \ge 5$ was inadvertently
omitted.

\cite{ZarhinMZ},  p. 538, Theorem 1.1, condition (ii): the $n$ should be $m$.

\end{document}